\documentclass[12pt,sort&compress]{article}
\usepackage{amssymb}
\usepackage{amsfonts}
\usepackage{amsmath}
\usepackage{amssymb,amsmath}
\usepackage{amscd}
\def\a{\alpha}
\def\g{\gamma}

\def\Vir{\hbox{Vir}}
\def\cl{\centerline}

\def\vs{\vspace*}
\def\W{\mathcal{HV}}

\def\Z{\mathbb{Z}}
\def\CC{\mathbb{C}}

\def\QED{\hfill$\Box$}

\def\t{\tilde}

\def\pa{\partial}
\def\la{\lambda}

\numberwithin{equation}{section}
\newtheorem{theo}{Theorem}[section]
\newtheorem{defi}[theo]{Definition}

\newtheorem{lemm}[theo]{Lemma}
\newtheorem{prop}[theo]{Proposition}

\newtheorem{rema}[theo]{Remark}

\makeatletter

\def\@biblabel#1{#1.~}

\makeatother

\begin{document}
\vs{10pt} \cl{\large {\bf Cohomology of Heisenberg-Virasoro conformal algebra}
\footnote{Corresponding author: Henan Wu (wuhenan@sxu.edu.cn).}} \vs{12pt}

\cl{ Lamei Yuan$^{\,\ddag}$, Henan Wu$^{\,\dag}$ }
 \cl{\small{ $^{\ddag}$ Academy of Fundamental and Interdisciplinary
 Sciences,}}\cl{\small{Harbin Institute of Technology, Harbin 150080, China}}
\cl{\small{$^{\dag}$School of Mathematical Sciences, Shanxi University, Taiyuan 030006, China}}
\cl{\small E-mail:
lmyuan@hit.edu.cn, wuhenan@sxu.edu.cn
 }\vs{6pt}
{\small\parskip .005 truein \baselineskip 3pt \lineskip 3pt

\noindent{\bf Abstract:} In this paper, we compute the cohomology of the Heisenberg-Virasoro conformal algebra with coefficients
in its modules, and in particular with trivial coefficients both for the basic and reduced complexes.
 \vs{5pt}

\vs{5pt}

\noindent{\bf MR(2000) Subject Classification:}~ 17B65, 17B68


\vs{18pt}

\noindent{\bf 1. \
Introduction}\setcounter{section}{1}\setcounter{theo}{0}\setcounter{equation}{0}
\vs{6pt}

The notion of Lie conformal algebra, introduced by Kac in \cite{K1}, encodes an axiomatic description
of the operator product expansion of chiral fields in conformal field
theory. In a more general context, a Lie conformal algebra is just an algebra in the pseudotensor category \cite{BDK}. Closely related to
vertex algebras, Lie conformal algebras have many applications in other areas of algebras and integrable systems. In particular, they give us powerful tools for the
study of infinite-dimensional Lie (super)algebras and associative algebras (and their representations),
satisfying the sole locality property \cite{K3}. The main examples of such Lie algebras
are those ¡°based¡± on the punctured complex plane, such as the Virasoro algebra and the loop Lie
algebras \cite{DK}. In addition, Lie conformal algebras resemble Lie algebras in many ways \cite{K2,SYX, SY,Z1,Z2}. A general cohomology theory of conformal algebras with coefficients in an arbitrary conformal module was
developed in \cite{BKV}, where explicit computations of cohomologies for the Virasoro conformal algebra and current conformal algebra were given.
 The low-dimensional cohomologies of the general Lie
conformal algebras $gc_N$ were studied in \cite{S}. The cohomologies of the $W(2,2)$-type conformal algebra with trivial coefficients were completely determined in \cite{YW}.

In this paper, we study the cohomology of the Heisenberg-Virasoro conformal algebra, which was introduced in \cite{SY} as a Lie conformal algebra associated with the twisted Heisenberg-Virasoro Lie algebra. By definition, the Heisenberg-Virasoro conformal algebra is free Lie conformal algebra of rank $2$ with a $\CC[\pa]$-basis $\{L,M\}$ and satisfying
\begin{eqnarray}\label{la-brac}
[L_\la L]=(\pa+2\la)L,\ \ [L_\la M]=(\pa+\la)M, \ \  [M_\la L]=\la M,\ \ [M_\la M]=0.
\end{eqnarray}
We denote by $\mathcal{HV}$ the Heisenberg-Virasoro conformal algebra. It is easily to see that $\mathcal{HV}$  contains the
Virasoro conformal algebra $\Vir$ as a subalgebra, which is a
free $\CC[\partial]$-module generated by $L$ such that
\begin{eqnarray}\label{Vir}
\Vir=\CC[\partial]L,\ \ \ [L_\lambda L]=(\partial+2\lambda)L.
\end{eqnarray}
Moreover, $\mathcal{HV}$  has a nontrivial
abelian conformal ideal with one free generator $M$ as a $\CC[\pa]$-module. Thus it is neither simple nor semi-simple.

The paper is organized as follows. In Section 2, we recall the notions of Lie conformal algebra, conformal module and cohomology of Lie conformal algebras, and present the main results of this paper (see Theorem \ref{main}). Section 3 is devoted to the proof of the main theorem.

Throughout this paper, all vector spaces and tensor products are over the complex field $\mathbb{C}$.  We use notations $\Z$ for the set of integers and $\Z_+$ for the set of nonnegative integers.
\vs{8pt}

\vs{8pt} \noindent{\bf 2. Preliminaries and Main results
}\setcounter{section}{2}\setcounter{theo}{0}\setcounter{equation}{0}\vs{6pt}

In this section, we recall the definition of a Lie conformal algebra and a (conformal) module over it and cohomology with coefficients in an arbitrary module. Then we list our main results of this paper.

\begin{defi} \rm (\cite{K1})
A Lie conformal algebra $\mathcal {R}$ is a $\CC[\partial ]$-module endowed with a $\CC$-bilinear map $\mathcal {R}\otimes \mathcal {R}\rightarrow \CC[\lambda]\otimes \mathcal {R}$, $a\otimes b \mapsto [a_\lambda b],
$ and
satisfying the following axioms ($a, b, c\in \mathcal {R}$),
\begin{eqnarray}
[\partial a_\lambda b]&=&-\lambda[a_\lambda b],\ \ [ a_\lambda \partial b]=(\partial+\lambda)[a_\lambda b] \ \ \mbox{(conformal\  sesquilinearity)},\label{Lc1}\\
{[a_\lambda b]} &=& -[b_{-\lambda-\partial}a] \ \ \mbox{(skew-symmetry)},\label{Lc2}\\
{[a_\lambda[b_\mu c]]}&=&[[a_\lambda b]_{\lambda+\mu
}c]+[b_\mu[a_\lambda c]]\ \ \mbox{(Jacobi \ identity)}\label{Lc3}.
\end{eqnarray}
\end{defi}
\begin{defi}\rm (\cite{CK}) A module $V$ over a Lie conformal algebra $\mathcal {A}$
is a $\mathbb{C}[\partial]$-module endowed with a $\CC$-bilinear map
$\mathcal {A}\otimes V\rightarrow
V[[\lambda]]$, $a\otimes v\mapsto a_\lambda v$, satisfying the following relations for $a,b\in\mathcal {A}$, $v\in V$,
\begin{eqnarray*}
&&a_\lambda(b_\mu v)-b_\mu(a_\lambda v)=[a_\lambda b]_{\lambda+\mu}v,\\
&&(\partial a)_\lambda v=-\lambda a_\lambda v,\ a_\lambda(\partial
v)=(\partial+\lambda)a_\lambda v.
\end{eqnarray*}
If $a_\lambda v\in
V[\lambda]$ for all $a \in \mathcal {A}$, $v \in V$, then $V$ is
called conformal. If $V$ is finitely generated over
$\mathbb{C}[\partial]$, then $V$ is simply called finite.
\end{defi}
Since we only consider conformal modules, we will simply shorten the term
``conformal module" to ``module". The vector space $\CC$ is viewed as a trivial module with trivial actions of $\pa$ and $\mathcal {A}$. For a fixed nonzero complex constant $a$, there
is a natural $\CC[\pa]$-module $\CC_a$, such that $\CC_a=\CC$ and $\pa v=a v$ for $v\in\CC_a$. Then
$\CC_a$ becomes an $\mathcal{A}$-module with
$\mathcal{A}$ acting by zero.

 For the Virasoro conformal algebra $\Vir$ (cf. \eqref{Vir}), it was shown in
\cite{CK} that all the free nontrivial $\Vir$-modules of rank $1$
over $\mathbb{C}[\partial]$ are the following ones $(\Delta,
\alpha\in \mathbb{C})$:
\begin{eqnarray}
M_{\Delta,\alpha}=\mathbb{C}[\partial]v,\ \ L_\lambda
v=(\partial+\alpha+\Delta \lambda)v.
\end{eqnarray}
The module $M_{\Delta,\alpha}$ is irreducible if and only if
$\Delta\neq 0$. The module $M_{0,\alpha}$ contains a unique
nontrivial submodule $(\partial +\alpha)M_{0,\alpha}$ isomorphic to
$M_{1,\alpha}.$  Moreover, the modules $M_{\Delta,\alpha}$ with
$\Delta\neq 0$ exhaust all finite irreducible nontrivial
$\Vir$-modules.

From the proof of \cite[Theorem 4.5 (1)]{SY}, we have
\begin{prop} \label{mod} All free nontrivial $\mathcal{HV}$-modules of rank 1 over
$\mathbb{C}[\pa]$ are the following ones:
\begin{eqnarray*}
M_{\Delta,\alpha,\beta}=\mathbb{C}[\partial]v,\ L_\lambda
v=(\partial+\alpha+\Delta \lambda)v, \ M_\lambda v=\beta v, \ \mbox{for
some}\ \Delta,\a,\beta\in\CC.
\end{eqnarray*}
\end{prop}

\begin{defi}\label{cochain}\rm (\cite{BKV}) An $n$-cochain ($n\in\Z_+$) of a Lie conformal algebra $\mathcal{A}$ with coefficients in an
$\mathcal{A}$-module $V$ is a $\CC$-linear map\vs{-5pt}
\begin{eqnarray*}
\gamma:\mathcal{A}^{\otimes n}\rightarrow V[\la_1,\cdots,\la_n],\ \
\ a_1\otimes\cdots \otimes a_n \mapsto
\g_{\la_1,\cdots,\la_n}(a_1,\cdots,a_n)
\end{eqnarray*}
satisfying the following conditions:\begin{itemize}\parskip-3pt
\item[\rm(1)] $\g_{\la_1,\cdots,\la_n}(a_1,\cdots,\pa a_i,\cdots,
a_n)=-\la_i\g_{\la_1,\cdots,\la_n}(a_1,\cdots, a_n)$ \ (conformal antilinearity),
\item[\rm (2)] $\g$ is skew-symmetric with respect to simultaneous permutations
of $a_i$'s and $\la_i$'s \ (skew-symmetry).
\end{itemize}
\end{defi}

As usual, let $\mathcal{A}^{\otimes 0}= \CC$, so that a $0$-cochain
is an element of $V$. Denote by  ${\t C}^n(\mathcal {A},V)$ the set
of all $n$-cochains. The differential $d$ of an $n$-cochain $\g$ is
defined as follows:
\begin{eqnarray}\label{ddd}
&&(d\g)_{\la_1,\cdots,\la_{n+1}}(a_1,\cdots,a_{n+1})\nonumber\\&&\ \ \ =\mbox{$\sum\limits_{i=1}^{n+1}$}(-1)^{i+1}a_{i_{\la_i}}\g_{\la_1,\cdots,\hat{\la_i},\cdots,\la_{n+1}}(a_1,\cdots,\hat{a_i},\cdots,a_{n+1})\nonumber\\
&&\ \ \ \ \ \ \ +\mbox{$\sum\limits_{i,j=1;i<j}^{n+1}$}(-1)^{i+j}\g_{\la_i+\la_j,\la_1,\cdots,\hat{\la_i},\cdots,\hat{\la_j},\cdots,\la_{n+1}}([a_{i_{\la_i}}a_j],a_1,\cdots,\hat{a_i},\cdots,\hat{a_j},\cdots,a_{n+1}),
\end{eqnarray}
where $\g$ is linearly extended over the polynomials in $\la_i$. In
particular, if $\g\in V$ is a $0$-cochain, then
$(d\g)_\la(a)=a_\la\g$.

It was shown in \cite{BKV} that the operator $d$ preserves the space of cochains and $d^2=0$. Thus the cochains of a Lie conformal algebra $\mathcal{A}$ with coefficients in an $\mathcal{A}$-module $V$ form a complex, called the {\it basic complex} and will be denoted by
\begin{eqnarray}
\t C^\bullet(\mathcal{A},V)=\mbox{$\bigoplus\limits_{
n\in\Z_+}$}\t C^n(\mathcal{A},V).
\end{eqnarray}
Moreover, define a (left) $\CC[\pa]$-module structure on $\t
C^\bullet(\mathcal{A},V)$ by
\begin{eqnarray*}
(\pa\g)_{\la_1,\cdots,\la_n}(a_1,\cdots,
a_n)=(\pa_V+\mbox{$\sum\limits_{i=1}^n$}\la_i)\g_{\la_1,\cdots,\la_n}(a_1,\cdots,
a_n),
\end{eqnarray*}
where $\pa_V$ denotes the action of $\pa$ on $V$. Then $d\pa=\pa d$
and thus $\pa \t C^\bullet(\mathcal{A},V)\subset \t
C^\bullet(\mathcal{A},V)$ forms a subcomplex. The quotient
complex
\begin{eqnarray*}
C^\bullet(\mathcal{A},V)=\t C^\bullet(\mathcal{A},V)/\pa \t
C^\bullet(\mathcal{A},V)= \mbox{$\bigoplus\limits_{n\in\Z_+}$}
C^n(\mathcal{A},V)
\end{eqnarray*}
is called the {\it reduced complex}.
\begin{defi}\label{def11}\rm The basic cohomology ${\rm \t H}^\bullet (\mathcal{A},V)$ of a Lie conformal algebra $\mathcal{A}$ with coefficients
 in an $\mathcal{A}$-module $V$ is the
cohomology of the basic complex $\t C^\bullet(\mathcal{A},V)$ and
the (reduced) cohomology ${\rm \,H}^\bullet (\mathcal{A},V)$ corresponds to
the reduced complex $C^\bullet(\mathcal{A},V)$.
\end{defi}

The following theorem is our main results of this paper.
\begin{theo}\label{main}  For the Heisenberg-Virasoro conformal algebra $\mathcal{HV}$, the
following statements hold.\begin{itemize}\parskip-3pt
\item[\rm(1)] For the trivial module $\CC$,
\begin{eqnarray}\label{h1}
&&{\rm dim\,\t H}^q(\mathcal{HV},\CC)=\left\{
\begin{array}{ll}
1 &{\mbox if}\ q=0,\\
3 &{\mbox if}\ q=3,\\
2 &{\mbox if}\ q=4,\\
0 &{\mbox otherwise},
\end{array}
\right.\\
&&
{\rm dim\, H}^q(\mathcal{HV},\CC)=\left\{
\begin{array}{ll}
1 &{\mbox if}\ q=0,\\
3 &{\mbox if}\ q=2,\\
5 &{\mbox if}\ q=3,\\
2 &{\mbox if}\ q=4,\\
0 &{\mbox otherwise}.
\end{array}
\right.\label{h2}
\end{eqnarray}
\item[\rm(2)] ${\rm H}^\bullet(\mathcal{HV},\CC_a)=0$ if $a\neq 0$.

\item[\rm(3)] ${\rm H}^\bullet(\mathcal{HV},M_{\Delta,\a})=0$ if $\a\neq 0$.
\end{itemize}
\end{theo}

\begin{rema} \rm Theorem \ref{main}(1) in particular shows that there is a unique nontrivial universal central extension of the Heisenberg-Virasoro conformal algebra $\mathcal{HV}$ by a three-dimensional center $\CC C_1\oplus\CC C_2\oplus \CC C_3$, which agrees with that of the twisted Heisenberg-Virasoro Lie algebra. The three independent reduced 2-cocycle $\bar \phi_1$, $\bar \phi_2$ and $\bar \phi_3$ are determined by \eqref{central-ex1}--\eqref{central-ex3} respectively, and the corresponding universal central extension $\widetilde{\mathcal{HV}}$ of $\mathcal{HV}$ is given by
\begin{eqnarray*}
&&[L_\la L]=(\pa+2\la)L+\frac{\la^3}{12} C_1,\\&&  [L_\la M]=(\pa+\la)M+\la^2 C_2, \\&& [M_\la L]=\la M-\la^2 C_2,\\&& [M_\la M]=\la C_3,
\end{eqnarray*}
where $C_1, C_2, C_3$ are nonzero central elements of $\widetilde{\mathcal{HV}}$ with $\pa C_i=0$, $i=1, 2, 3$.
\end{rema}

\begin{rema}\rm  Denote by ${\rm Lie}(\W)_-$ the annihilation Lie algebra of $\W$.
It can be easily checked that ${\rm Lie}(\W)_-$ is isomorphic to the subalgebra spanned by
$\{L_n, M_n\big| -1\leq n\in\Z\}$ of the centerless Heisenberg-Virasoro algebra.
Since ${\rm \t H}^q(\mathcal{HV},\CC)\cong{\rm H}^q({\rm Lie}(\mathcal{HV})_-,\CC)$, we have actually determined
the cohomology group of ${\rm Lie}(\mathcal{HV})_-$ with trivial coefficients (cf. \cite{BKV}).
\end{rema}

\vs{8pt} \noindent{\bf 3. Proof of Theorem \ref{main}
}\setcounter{section}{3}\setcounter{theo}{0}\setcounter{equation}{0}
\vs{6pt}

In this section, we prove Theorem \ref{main}, which will be done by several lemmas.

Keep notations in the previous section. For $\gamma\in{\t C}^q(\mathcal {A},V)$, we call
$\gamma$ a {\it $q$-cocycle} if $d(\gamma)=0$; a {\it $q$-coboundary} if there exists a $(q-1)$-cochain $\phi\in\t
C^{q-1}(\mathcal{A},V)$ such that $\gamma=d(\phi)$. Two cochains
$\gamma_1$ and $\gamma_2$ are called {\it equivalent} if $\gamma_1-\gamma_2$ is a coboundary. Denote by $\t D^q(\mathcal{A},V)$ and $\t B^q(\mathcal{A},V)$ the spaces
of $q$-cocycles and $q$-boundaries, respectively. By Definition \ref{def11},
\begin{eqnarray*}
{\rm \t H}^q(\mathcal{A},V)=\t D^q(\mathcal{A},V)/\t B^q(\mathcal{A},V)=\{\mbox{equivalent classes of
$q$-cocycles}\}.
\end{eqnarray*}

\begin{lemm} ${\rm \t
H}^0(\mathcal{HV},\CC)={\rm H}^0(\mathcal{HV},\CC)=\CC$.
\end{lemm}
\noindent{\it Proof.~} For any $\gamma\in \t
C^0(\mathcal{HV},\CC)=\CC$, $(d\gamma)_\la (X)=X_\la \gamma =0$ for
$X\in \mathcal{HV}$. This means $\t D^0(\mathcal{HV},\CC)=\CC$ and $\t B^0(\mathcal{HV},\CC)=0$. Thus ${\rm \t
H}^0(\mathcal{HV},\CC)=\CC$ and $ {\rm H}^0(\mathcal{HV},\CC)=\CC$
since $\pa\CC=0$.\QED\vskip5pt

 Let $\gamma\in \t
C^q(\mathcal{HV},\CC)$ with $q>0$. By Definition \ref{cochain}, $\g$ is determined by its value on $X_1\otimes\cdots \otimes X_q$ with $X_i\in\{L,M\}$. Without loss of generality, we always assume that the first $k$ variables are $L$ and the last $q-k$ variables are $M$ in $\g_{\la_1,\cdots,\la_{q}}(X_1,\cdots,X_q)$, since $\g$ is skew-symmetric.
Thus we can regard $\g_{\la_1,\cdots,\la_{q}}(X_1,\cdots,X_q)$ as a polynomial in $\la_1,\cdots,\la_q$, which is skew-symmetric in $\la_1,\cdots,\la_k$ and also skew-symmetric in $\la_{k+1},\cdots,\la_q$. Therefore, $\g_{\la_1,\cdots,\la_{q}}(X_1,\cdots,X_q)$ is divisible by $$\mbox{$\prod\limits_{1\leq i< j\leq k}$}(\la_i-\la_j)\times\mbox{$\prod\limits_{k+1\leq i< j\leq q}$}(\la_i-\la_j),$$ whose polynomial degree is $k(k-1)/2+(q-k)(q-k-1)/2$.

Following \cite{BKV}, we define an operator $\tau:\t C^q(\mathcal{HV},\CC)\rightarrow
\t C^{q-1}(\mathcal{HV},\CC)$ by
\begin{eqnarray}\label{7++}
(\tau
\g)_{\la_1,\cdots,\la_{q-1}}(X_1,\cdots,X_{q-1})=(-1)^{q-1}\frac{\pa}{\pa\la}\g_{\la_1,\cdots,\la_{q-1},\la}(X_1,\cdots,X_{q-1},L)|_{\la=0},
\end{eqnarray}
where $X_1=\cdots=X_{k}=L,$ $ X_{k+1}=\cdots =X_{q-1}=M$.  By \eqref{ddd}, \eqref{7++} and skew-symmetry of $\g$,
\begin{eqnarray}\label{7++1}
&&((d\tau+\tau d)
\g)_{\la_1,\cdots,\la_{q}}(X_1,\cdots,X_q)\nonumber\\
&&\ \ \ \ =(-1)^q\frac{\pa}{\pa\la}\mbox{$\sum\limits_{i=1}^q$}(-1)^{i+q+1}\g_{\la_i+\la,\la_1,\cdots,\hat{\la_i},\cdots,\la_{q}}([X_{i\,{\la_i}} L],X_1,\cdots,\hat{X_{i}},\cdots, X_q)|_{\la=0}\nonumber\\
&&\ \ \ \ =\frac{\pa}{\pa\la}\mbox{$\sum\limits_{i=1}^q$}\g_{\la_1,\cdots,\la_{i-1},\la_i+\la,\la_{i+1},\cdots,\la_{q}}
(X_1,\cdots,X_{i-1},[X_{i\,{\la_i}} L],X_{i+1},\cdots,X_q)|_{\la=0}.
\end{eqnarray}
By \eqref{la-brac} and conformal antilinearity of $\g$, $[X_{i\,{\la_i}} L]$ can be replaced by either
$(\la_i-\la)X_i$ when $X_i=L$ or by
$\la_i X_i$ when $X_i=M$ in \eqref{7++1}.
Thus, equality \eqref{7++1} becomes
\begin{eqnarray}\label{7++2}
&&((d\tau+\tau d)
\g)_{\la_1,\cdots,\la_{q}}(X_1,\cdots,X_q)\nonumber\\
&&\ \ \ \ =\frac{\pa}{\pa\la}\mbox{$\sum\limits_{i=1}^k$}(\la_i-\la) \g_{\la_1,\cdots,\la_{i-1},\la_i+\la,\la_{i+1},\cdots,\la_{q}}
(X_1,\cdots,X_{i-1},X_i,X_{i+1},\cdots,X_q)|_{\la=0}\nonumber\\
&&\ \ \ \ \ \ \,+ \frac{\pa}{\pa\la}\mbox{$\sum\limits_{i=k+1}^q$}\la_i \g_{\la_1,\cdots,\la_{i-1},\la_i+\la,\la_{i+1},\cdots,\la_{q}}
(X_1,\cdots,X_{i-1},X_i,X_{i+1},\cdots,X_q)|_{\la=0}\nonumber\\
&&\ \ \ \ =({\rm deg\,} \g-k)\g_{\la_1,\cdots,\la_{q}}(X_1,\cdots,X_q),
\end{eqnarray}
where ${\rm deg\,}\g$ is the total degree of $\g$ in $\la_1,\cdots,\la_q$. Therefore, only those homogeneous cochains whose degree as a polynomial is equal to $k$ contribute to the cohomology of $\t C^\bullet(\mathcal{HV},\CC)$. Consider the quadratic inequality $$\frac{k(k-1)}{2}+\frac{(q-k)(q-k-1)}{2}\leq k,$$ whose discriminant is $\Delta_k=-4k^2+12k+1$. Since $\Delta_k\geq 0$ has $k=0, 1, 2$ and $3$ as the only integral solutions, we have
\begin{eqnarray}\label{qqq}
k=\left\{
\begin{array}{ll}
0, 1  & if \ q=1,\\
1, 2  & if \ q=2,\\
1, 2, 3  & if \ q=3, \\
3, 4  & if \ q=4.
\end{array}
\right.
\end{eqnarray}
In particular, ${\rm \t H}^q(\mathcal{W},\CC)=0$ for $q\geq 5$.

\begin{lemm} Theorem \ref{main} (1) holds.
\end{lemm}
\noindent{\it Proof.~} It needs to compute ${\rm \t H}^q(\mathcal{HV},\CC)$ for $q=1,2,3,4$.

For $q=1$, we need to consider $k=0, 1$ by \eqref{qqq}. Let $\g$ be a 1-cocycle.
From the discussions below \eqref{7++2}, we know $\g_\la(M)$ should be a constant, whereas $\g_\la(L)$ should be a constant factor of $\la$. By $d\g=0$, it is easy to check that both $\g_\la(M)$ and $\g_\la(L)$ are zero.
Hence,  ${\rm \t
H}^1(\mathcal{HV},\CC)=0$.

For $q=2$, we need to consider $k=1, 2$ by \eqref{qqq}. If $\g\in {\rm \t
H}^2(\mathcal{HV},\CC)$, then ${\rm deg}(\g_{\la_1,\la_2}(L,L))=2$ and ${\rm deg}(\g_{\la_1,\la_2}(L,M))=1$ as polynomials in $\la_1,\la_2$. By skew-symmetry of $\g$, $\g_{\la_1,\la_2}(L,L)$ should be a constant factor of $\la_1^2-\la_2^2$, which is a coboundary of a $1$-cochain of the form $\varphi_{\la_1}(L)=\la_1$.  Assume that $\g_{\la_1,\la_2}(L,M)=a\la_1+b\la_2$ for some $a,b\in\CC$.
A straightforward computation shows that
\begin{equation}\label{hq2}
(d\g)_{\la_1,\la_2,\la_3}(L,L,M)=-a\la_1(\la_1+\la_2)(\la_1-\la_2)+a\la_2\la_3-a\la_1\la_3.
\end{equation}
Then $d\g=0$ gives $a=0$. Set $\varphi_{\la_1}(M)=b$. Then $(d\varphi)_{\la_1,\la_2}(L,M)=b\la_2=\g_{\la_1,\la_2}(L,M)$, namely, $\g_{\la_1,\la_2}(L,M)$ is also a coboundary. Thus ${\rm \t
H}^2(\mathcal{HV},\CC)=0$.

 For $q=3$, we need to consider $k=1,2,3$ by \eqref{qqq}. Let $\g\in\t D^3(\W,\CC)$ be a 3-cocycle. In the case of $k=1$,  we can assume that $\gamma_{\la_1,\la_2,\la_3}(L,M,M)=c(\la_2-\la_3)$ for some $c\in\CC$. One can check that it satisfies the following equation
\begin{eqnarray}\label{q=3-4}
&&(d\gamma)_{\la_1,\la_2,\la_3,\la_4}(L,L,M,M)\nonumber\\=&&c(-(\la_1-\la_2)\phi_{\la_1+\la_2,\la_3,\la_4}(L,M,M)+\la_3\phi_{\la_2,\la_1+\la_3,\la_4}(L,M,M)\nonumber\\
&&-\la_4\phi_{\la_2,\la_1+\la_4,\la_3}(L,M,M)-\la_3\phi_{\la_1,\la_2+\la_3,\la_4}(L,M,M)+\la_4\phi_{\la_1,\la_2+\la_4,\la_3}(L,M,M))\nonumber\\
=&&c(-(\la_1-\la_2)(\la_3-\la_4)+\la_3(\la_1+\la_3-\la_4)\nonumber\\&&-\la_4(\la_1+\la_4-\la_3)-\la_3(\la_2+\la_3-\la_4)+\la_4(\la_2+\la_4-\la_3))\nonumber\\
=&&0.
\end{eqnarray}
And it is a not coboundary, because it can be the coboundary of a 2-cochain $\varphi_{\la_1,\la_2}(M,M)$ of degree 0, which must be zero by skew-symmetry of $\varphi$. In the case when $k=2$, we suppose that
\begin{eqnarray}\label{q=3-5}
\gamma_{\la_1,\la_2,\la_3}(L,L,M)=(\la_1-\la_2)(a(\la_1+\la_2)+b\la_3), \ {\mbox {for some}}\ a,b\in\CC.
\end{eqnarray}
It satisfies $(d\gamma)_{\la_1,\la_2,\la_3,\la_4}(L,L,L,M)=0.$
Setting $\varphi_{\la_1,\la_2}(L,M)=a\la_1$, we have $$(d\varphi)_{\la_1,\la_2,\la_3}(L,L,M)=-a(\la_1-\la_2)(\la_1+\la_2+\la_3),$$
and
\begin{eqnarray*}
(d\varphi)_{\la_1,\la_2,\la_3}(L,L,M)+\gamma_{\la_1,\la_2,\la_3}(L,L,M)=(b-a)(\la_1-\la_2)\la_3.
\end{eqnarray*}
 Thus $\gamma_{\la_1,\la_2,\la_3}(L,L,M)$ is equivalent to $(\la_1-\la_2)\la_3$, which is not a coboundary by (\ref{hq2}). According to \cite[Theorem 7.1]{BKV},
$\gamma_{\la_1,\la_2,\la_3}(L,L,L)=(\la_1-\la_2)(\la_1-\la_3)(\la_2-\la_3)$ (up to a constant factor) is also a 3-cocycle, which is not a coboundary. Therefore, ${\rm dim\, H}^3(\W,\CC)=3$ and ${\rm \t H}^3(\W,\CC)=\CC \phi_1 \oplus \CC \phi_2 \oplus \CC \phi_3$, where
\begin{eqnarray}
&&\phi_1:=\phi_{1\,{\la_1,\la_2,\la_3}}(L,M,M)=\la_2-\la_3,\label{hvq-1}\\
&&\phi_2:=\phi_{2\,{\la_1,\la_2,\la_3}}(L,L,M)=(\la_1-\la_2)\la_3,\label{hvq-2}\\
&&\phi_3:=\phi_{3\,{\la_1,\la_2,\la_3}}(L,L,L)=(\la_1-\la_2)(\la_1-\la_3)(\la_2-\la_3),\label{hvq-3}
\end{eqnarray}
and where, we take $\phi_1$ for example, the skew-symmetric function $\phi_1:\W\otimes\W\otimes\W\rightarrow \CC[\la_1,\la_2,\la_3]$ has values $\la_2-\la_3$ on $L\otimes M\otimes M$ and $0$ on others.

For $q=4$, we need to consider $k=2, 3$. Let  $\g\in\t D^4(\W,\CC)$. By skew-symmetry of $\g$ with whose degree as a polynomial taken into account, we assume that
 \begin{eqnarray}
  \g_{\la_1,\la_2,\la_3,\la_4}(L,L,M,M)=e_1(\la_1-\la_2)(\la_3-\la_4), \ \ for\ some\ \  e_1\in\CC.\label{q=4-1}
 \end{eqnarray}
We have $(d\g)_{\la_1,\la_2,\la_3,\la_4,\la_5}(L,L,L,M,M)=0.$
 Thus $\psi_1:=\psi_{1\,{\la_1,\la_2,\la_3,\la_4}}(L,L,M,M)=(\la_1-\la_2)(\la_3-\la_4)$ is a $4$-cocycle. If $\psi_1$ is the coboundary of a 3-cochain $\phi_{\la_1,\la_2,\la_3}(L,M,M)$ of degree 1, then $\phi_{\la_1,\la_2,\la_3}(L,M,M)$ must be a constant factor of $\la_2-\la_3$, whose coboundary is zero by \eqref{q=3-4}.
Hence $\psi_1$ is not a coboundary. Similarly, assume that
 \begin{eqnarray}
  \g_{\la_1,\la_2,\la_3,\la_4}(L,L,L,M)=e_2(\la_1-\la_2)(\la_2-\la_3)(\la_1-\la_3),\ \ for \ some \ e_2\in\CC.\label{q=4-2}
 \end{eqnarray}
One can check that it satisfies $d\g=0$. We should point out that $\psi_2:=\psi_{2\,{\la_1,\la_2,\la_3,\la_4}}(L,L,M,M)=(\la_1-\la_2)(\la_2-\la_3)(\la_1-\la_3)$ is not a coboundary. Because it can be the coboundary of a 3-cochain $\phi_{\la_1,\la_2,\la_3}(L,L,M)$ of degree 2, which should be of form $(\la_1-\la_2)(a(\la_1+\la_2)+b\la_3)$, whose coboundary is zero by \eqref{q=3-5}.
 Hence ${\rm \tilde H}^4(\W,\CC)=\CC \psi_1\oplus \CC \psi_2$.

 By \cite[Proposition 2.1]{BKV}, the map $\gamma\mapsto\partial \gamma$ gives an isomorphism ${\rm \t H}^q(\W,\CC)\cong {\rm H}^q(\pa\t C^\bullet) $ for $q\geq 1$.
Therefore,
\begin{eqnarray}\label{pa-c}
{\rm H}^q(\pa\t C^\bullet)=\left\{
\begin{array}{ll} \CC (\pa\phi_1)\oplus\CC (\pa\phi_2)\oplus\CC (\pa\phi_3) &{if}\ q=3,\\
\CC (\pa\psi_1)\oplus \CC(\pa\psi_2) &{if}\ q=4,\\
0 &{otherwise}.
\end{array}
\right.
\end{eqnarray}
The computation of ${\rm H}^\bullet(\W,\CC)$ is based on
the short exact sequence of complexes
 \begin{equation}\label{exact}
\begin{CD}
0@>>> {\pa\t C^\bullet} @>{\rm \iota}>> {\t C^\bullet} @>{\rm \pi}>> C^\bullet @>>> 0
\end{CD}
\end{equation}
 where $\iota$  and $\rm \pi$ are the embedding and the natural projection, respectively. The exact sequence \eqref{exact} gives the following long exact sequence of cohomology groups (cf. \cite{BKV}):
\begin{equation}\label{key}
\begin{CD}
\cdots @>>> {\rm H}^q(\pa\t C^\bullet) @>{\rm \iota}_q>> {\rm \t H}^q(\W,\CC)  @>{\rm \pi}_q>> {\rm H}^q(\W,\CC) @>{\rm \omega}_q>>\\
@>>> {\rm H}^{q+1}(\pa\t C^\bullet) @>{\rm \iota}_{q+1}>> {\rm \t H}^{q+1}(\W,\CC)  @>{\rm \pi}_{q+1}>> {\rm H}^{q+1}(\W,\CC) @>>>\cdots \\
\end{CD}
\end{equation}
where $\iota_q, \pi_{q}$ are induced by $\iota, \pi$ respectively and $w_q$ is the $q-$th connecting hommorphism.
Given $\partial \gamma\in {\rm H}^q(\pa\t C^\bullet)$ with a nonzero element $\gamma\in {\rm \t H}^q(\W,\CC)$ of degree $k$,
we have $\iota_q(\partial \gamma)=\partial \gamma\in {\rm \t H}^q(\W,\CC)$.
Since ${\rm deg\,}(\partial\gamma)={\rm deg\,}(\gamma)+1=k+1$, $\partial\gamma=0\in {\rm \t H}^q(\W,\CC)$.
Thus the image ${\rm im}({\rm \iota}_{q})$ of $\iota_q$ is zero for any $q\in\Z_+$.
Because ${\rm ker}({\rm \pi}_{q})={\rm im}({\rm \iota}_{q})=\{0\}$ and ${\rm im}({\rm \omega}_{q})={\rm ker}({\rm \iota}_{q+1})={\rm H}^{q+1}(\pa\t C^\bullet)$,
we obtain the following short exact sequence
\begin{equation}\label{key1}
\begin{CD}
0 @>>> {\rm \t H}^q(\W,\CC)  @>{\rm \pi}_{q}>> {\rm H}^q(\W,\CC) @>{\rm \omega}_q>>{\rm H}^{q+1}(\pa\t C^\bullet)  @>>>0.
\end{CD}
\end{equation}
Therefore,
\begin{eqnarray}\label{key111}
{\rm dim\, H}^q(\W,\CC)={\rm dim\,\t H}^q(\W,\CC)+{\rm dim\, H}^{q+1}(\pa\t C^\bullet), \ \mbox{for\ all} \ q\geq 0.
\end{eqnarray}
Consequently, we obtain \eqref{h2}.
By \eqref{key1},
the basis of ${\rm H}^q(\W,\CC)$ can be obtained by combining the images of a basis of ${\rm \t H}^q(\W,\CC)$
with the pre-images of a basis of ${\rm \t H}^{q+1}(\W,\CC)$.
Let $\varphi$ be a nonzero $(q+1)$-cocycle of degree $k$ such that
$\partial\varphi\in{\rm H}^{q+1}(\pa\t C^\bullet)$.
By \eqref{7++2}, we have
\begin{eqnarray}
d(\tau(\pa\varphi)=(d\tau+\tau d)(\pa\varphi)=({\rm deg\,} (\pa\varphi)-k)(\pa\varphi)=((k+1)-k)(\pa\varphi)=\pa\varphi.
\end{eqnarray}
Thus the pre-image $\omega_q^{-1}(\partial\varphi)$ of $\pa\varphi$ under the connecting homorphism $\omega_p$ is $\tau(\pa\varphi)$, i.e., $\omega_q^{-1}(\partial\varphi)=\tau(\pa\varphi)$. Finally, let us finish the proof by giving the basis of ${\rm H}^q(\W,\CC)$ for $q=2,3,4$.
For $q=2$, we have known that ${\rm \t H}^2(\W,\CC)=0$ and ${\rm H}^3(\pa\t C^\bullet)= \CC (\pa\phi_1)\oplus\CC (\pa\phi_2)\oplus\CC (\pa\phi_3)$. By \eqref{7++}, \eqref{hvq-1}--\eqref{hvq-3},
\begin{eqnarray}
\bar\phi_1:&=&(\tau(\pa\phi_1))_{\la_1,\la_2}(M,M)\nonumber\\
&=&(-1)^2\frac{\pa}{\pa\la}(\pa\phi_1))_{\la_1,\la_2,\la}(M,M,L)|_{\la=0}\nonumber\\
&=&\frac{\pa}{\pa\la}(\la_1+\la_2+\la)(\la_1-\la_2)|_{\la=0}\nonumber\\
&=&\la_1-\la_2,\label{central-ex1}\\
\bar\phi_2:&=&(\tau(\pa\phi_1))_{\la_1,\la_2}(L,M)\nonumber\\
&=&(-1)^2\frac{\pa}{\pa\la}(\pa\phi_1))_{\la_1,\la_2,\la}(L,M,L)|_{\la=0}\nonumber\\
&=&-\frac{\pa}{\pa\la}(\la_1+\la_2+\la)(\la_1-\la)\la_2|_{\la=0}\nonumber\\
&=&\la_2^2,\label{central-ex2}
\end{eqnarray}
\begin{eqnarray}
\bar\phi_3:&=&(\tau(\pa\Lambda_3))_{\la_1,\la_2}(L,L)\nonumber\\
&=&(-1)^2\frac{\pa}{\pa\la}(\pa\Lambda_3))_{\la_1,\la_2,\la}(L,L,L)|_{\la=0}\nonumber\\
&=&\frac{\pa}{\pa\la}(\la_1+\la_2+\la)(\la_1-\la_2)(\la_2-\la)(\la_1-\la)|_{\la=0}\nonumber\\
&=&-\la_1^3+\la_2^3.\label{central-ex3}
\end{eqnarray}
Thus ${\rm H}^2(\W,\CC)=\CC \bar \phi_1\oplus \CC \bar \phi_2 \oplus \CC \bar \phi_3.$ For $q=3$, we have
\begin{eqnarray*}
\bar \psi_1:&=&(\tau(\pa\psi_1))_{\la_1,\la_2,\la_3}(L,M,M)\\&=&(-1)^3\frac{\pa}{\pa\la}(\pa\psi_1))_{\la_1,\la_2,\la_3,\la}(L,M,M,L)|_{\la=0}\\
&=&\frac{\pa}{\pa\la}(\la_1+\la_2+\la_3+\la)(\la_1-\la)(\la_3-\la_2)|_{\la=0}\\
&=&\la_2^2-\la_3^2,\\
\bar \psi_2:&=&(\tau(\pa\psi_2))_{\la_1,\la_2,\la_3}(L,L,M)\\&=&(-1)^3\frac{\pa}{\pa\la}(\pa\psi_2))_{\la_1,\la_2,\la_3,\la}(L,L,M,L)|_{\la=0}\\
&=&\frac{\pa}{\pa\la}(\la_1+\la_2+\la_3+\la)(\la_1-\la_2)(\la_1-\la)(\la_2-\la)|_{\la=0}\\
&=&-\la_1^3-\la_1^2\la_3+\la_2^3+\la_2^2\la_3.
\end{eqnarray*}
Hence, ${\rm H}^3(\W,\CC)=\CC \phi_1\oplus\CC \phi_2\oplus\CC\phi_3 \oplus \CC \bar\psi_1  \oplus \CC \bar\psi_2$ and ${\rm H}^4(\W,\CC)=\CC \psi_1\oplus \CC \psi_2$ by \eqref{key111}. \QED
\begin{lemm} Theorem \ref{main} (2) holds.
\end{lemm}
\noindent{\it Proof.~} For $a\neq 0$, define an operator $\tau_2:\t C^q(\W,\CC_a)\rightarrow
\t C^{q-1}(\W,\CC_a)$ by
\begin{eqnarray}\label{7-3}
(\tau_2
\g)_{\la_1,\cdots,\la_{q-1}}(X_1,\cdots,X_{q-1})=(-1)^{q-1}\g_{\la_1,\cdots,\la_{q-1},\la}(X_1,\cdots,X_{q-1},L)|_{\la=0},
\end{eqnarray}
for $X_1,\cdots,X_{q-1}\in\{L,M\}$. By the fact that $\pa\t C^q(\W,\CC_a)=(a+\sum_{i=1}^q\la_i)\t
C^q(\W,\CC_a)$, we
have
\begin{eqnarray}\label{7-4}
((d\tau_2+\tau_2 d)
\g)_{\la_1,\cdots,\la_{q}}(X_1,\cdots,X_q)&=&\big(\mbox{$\sum_{i=1}^q$}\la_i\big)\g_{\la_1,\cdots,\la_{q}}(X_1,\cdots,X_q)\nonumber\\
&\equiv& -a \g_{\la_1,\cdots,\la_{q}}(X_1,\cdots,X_q) \ (\mbox{mod}\
\pa\t C^q(\W, \CC_a).
\end{eqnarray}
Let $\g\in\t C^q(\mathcal{W},\CC_a)$ be
a $q$-cochain such that $d\g\in \pa\t C^{q+1}(\W,\CC_a)$, namely,
there is a $(q+1)$-cochain $\phi$ such that
$d\g=(a+\sum_{i=1}^{q+1}\la_i)\phi$. By \eqref{7-3}, we have $\tau_2 d\g=(a+\sum_{i=1}^q\la_i)\tau_2\phi\in\pa\t
C^q(\W,\CC_a)$. It follows from \eqref{7-4} that
$\g\equiv-d(a^{-1}\tau_2\g)$ is a reduced coboundary.\QED
\begin{lemm}
Theorem \ref{main}(3) holds.
\end{lemm}
\noindent{\it Proof.~} In this case, $\pa\t
C^q(\W,M_{\Delta,\a})=(\pa+\sum_{i=1}^q\la_i)\t C^q(\W,
M_{\Delta,\a})$. Define an
operator $\tau_3:C^q(\W,M_{\Delta,\a})\rightarrow
C^{q-1}(\W,M_{\Delta,\a})$ by
\begin{eqnarray*}\label{5-3}
(\tau_3
\g)_{\la_1,\cdots,\la_{q-1}}(X_1,\cdots,X_{q-1})=(-1)^{q-1}\g_{\la_1,\cdots,\la_{q-1},\la}(X_1,\cdots,X_{q-1},L)|_{\la=0},
\end{eqnarray*}
for $X_1,\cdots,X_{q-1}\in\{L,M\}$. We have
\begin{eqnarray}\label{5-4}
((d\tau_3+\tau_3 d)
\g)_{\la_1,\cdots,\la_{q}}(X_1,\cdots,X_q)&=&L_\la\g_{\la_1,\cdots,\la_{q}}(X_1,\cdots,X_{q})|_{\la=0}+
\big(\mbox{$\sum\limits_{i=1}^q$}\la_i\big)\g_{\la_1,\cdots,\la_{q}}(X_1,\cdots,X_q)\nonumber\\[-6pt]
&=&\big(\pa+\a+\mbox{$\sum\limits_{i=1}^q$}\la_i\big)\g_{\la_1,\cdots,\la_{q}}(X_1,\cdots,X_q)\nonumber\\[-6pt]
&\equiv& \a \g_{\la_1,\cdots,\la_{q}}(X_1,\cdots,X_q) \ (\mbox{mod}\
\pa\t C^q(\W, M_{\Delta,\a})).
\end{eqnarray}
If $\g$ is a reduced $q$-cocycle, then, by \eqref{5-4}, $\g\equiv d(\a^{-1}\tau_3\g)$ is a reduced coboundary, since $\a\neq
0$. Thus $ {\rm H}^q(\W,M_{\Delta,\a})=0$ for all $q\geq 0$.\QED\vskip5pt
This completes the proof of Theorem \ref{main}.

\vspace{4mm} \noindent\bf{\footnotesize Acknowledgements.}\ \rm
{\footnotesize This work was supported by National Natural Science
Foundation grants of China (11301109, 11526125) and the Research Fund for the Doctoral Program of Higher Education (20132302120042).}\\
\vskip18pt \small\footnotesize
\parskip0pt\lineskip1pt

\end{document}